\documentclass{amsart}%
\usepackage{amsfonts}
\usepackage{amsmath}
\usepackage{amssymb}
\usepackage{graphicx}%
\providecommand{\U}[1]{\protect\rule{.1in}{.1in}}
\newtheorem{theorem}{Theorem}
\theoremstyle{plain}

\newtheorem{definition}{Definition}

\newtheorem{lemma}{Lemma}

\newtheorem{remark}{Remark}

\numberwithin{equation}{section}
\begin{document}
\title{A Nonlinear Extension of Korovkin's Theorem}
\author{Sorin G. Gal}
\address{Department of Mathematics and Computer Science\\
University of Oradea\\
University\ Street No. 1, Oradea, 410087, Romania and Academy of Romanian Scientists, Splaiul Independentei No. 54, Bucharest
050094, Romania, e-mail : galso@uoradea.ro}
\author{Constantin P. Niculescu}
\address{Department of Mathematics, University of Craiova\\
Craiova 200585, Romania}
\curraddr{Academy of Romanian Scientists, Splaiul Independentei No. 54, Bucharest
050094, Romania, e-mail : cpniculescu@gmail.com}
\date{Preliminary version December 2,  2019}
\dedicatory{Dedicated to Professor Nicolae Dinculeanu, on the occasion of his 95th birthday.}
\begin{abstract}
In this paper we extend the classical Korovkin theorems to the framework of
comonotone additive, sublinear and monotone operators. Based on the theory of
Choquet capacities, several concrete examples illustrating our results are
also discussed.

\end{abstract}
\maketitle

\section{Introduction}

One of the most elegant results in the theory of approximation is Korovkin's
theorem, that provides a generalization of the well-known proof of
Weierstrass's classical approximation theorem as was given by Bernstein.

\begin{theorem}
\label{Thm1}\emph{(Korovkin} \emph{\cite{Ko1953}}, \emph{\cite{Ko1960})} Let
$(L_{n})_{n}$ be a sequence of positive linear operators that map $C([0,1])$
into itself. Suppose that the sequence $(L_{n}(f))_{n}$ converges to $f$
uniformly on $[0,1]$ for each of the test functions $1,$ $x$ and $x^{2}$. Then
this sequence converges to $f$ uniformly on $[0,1]$ for every $f\in C([0,1])$.
\end{theorem}

Simple examples show that the assumption concerning the positivity of the
operators $L_{n}$ cannot be dropped. What about the assumption on their linearity?

Over the years, many generalizations of Theorem 1 appeared, in a variety of
settings including important Banach function spaces. A nice account on the
present state of art is offered by the authoritative monograph of Altomare and
Campiti \cite{AC1994} and the excellent survey of Altomare \cite{Alt2010}. The
literature concerning the subject of Korovkin type theorems is really huge, a
search on Google offering more than 26,000 results. However, except for
Theorem 2.7 in the 1973 paper of Bauer \cite{Bauer}, the extension of this
theory beyond the framework of linear functional analysis remained largely unexplored.

Inspired by the Choquet theory of integrability with respect to a nonadditive
measure, we will prove in this paper that the restriction to the class of
positive linear operators can be relaxed by considering operators that verify
a mix of conditions characteristic for Choquet's integral.

As usually, for $X$ a Hausdorff topological space we will denote by
$\mathcal{F}(X)$ the vector lattice of all real-valued functions defined on
$X$ endowed with the pointwise ordering. Two important vector sublattices of
it are
\[
C(X)=\left\{  f\in\mathcal{F}(X):\text{ }f\text{ continuous}\right\}  \text{ }%
\]
and
\[
C_{b}(X)=\left\{  f\in\mathcal{F}(X):\text{ }f\text{ continuous and
bounded}\right\}  .
\]
With respect to the the sup norm, $C_{b}(X)$ becomes a Banach lattice. See
\cite{Sch1974} for the theory of these spaces.

Suppose that $X$ and $Y$ are two Hausdorff topological spaces and $E$ and $F$
are respectively vector sublattices of $C(X)$ and $C(Y).$ An operator
$T:E\rightarrow F$ is called:

- \emph{sublinear} if it is both \emph{subadditive}, that is,%
\[
T(f+g)\leq T(f)+T(g)\quad\text{for all }f,g\in E,
\]
and \emph{positively homogeneous,} that is,
\[
T(af)=aT(f)\quad\text{for all }a\geq0\text{ and }f\in E;
\]

- \emph{monotonic} if $f\leq g$ in $E$ implies $T(f)\leq T(g);$

- \emph{comonotonic additive} if $T(f+g)=T(f)+T(g)$ whenever the functions
$f,g\in E$ are comonotone in the sense that
\[
(f(s)-f(t))\cdot(g(s)-g(t))\geq0\text{ for all }s,t\in X.
\]

Our main result extends Korovkin's results to the framework of operators
acting on vector lattices of functions of several variables that play the
properties of sublinearity, monotonicity and comonotonic additivity. We use
families of test functions including the canonical projections on
$\mathbb{R}^{N}$,
\[
\operatorname*{pr}\nolimits_{k}:(x_{1},...,x_{N})\rightarrow x_{k}%
,\text{\quad}k=1,...,N.
\]

\begin{theorem}
\label{Thm2} \emph{(The nonlinear extension of Korovkin's theorem: the several
variables case)} Suppose that $X$ is a locally compact subset of the Euclidean
space $\mathbb{R}^{N}$ and $E$ is a vector sublattice of $\mathcal{F}(X)$ that
contains the test functions $1,~\pm\operatorname*{pr}_{1},...,~\pm
\operatorname*{pr}_{N}$ and $\sum_{k=1}^{N}\operatorname*{pr}_{k}^{2}$.

$(i)$ If $(T_{n})_{n}$ is a sequence of monotone and sublinear operators from
$E$ into $E$ such that
\begin{equation}
\lim_{n\rightarrow\infty}T_{n}(f)=f\text{\quad uniformly on the compact
subsets of }X \label{limcond}%
\end{equation}
for each of the $2N+2$ aforementioned test functions, then the property
\emph{(\ref{limcond})} also holds for all nonnegative functions $f$ in $E\cap
C_{b}(X)$.

$(ii)$ If, in addition, each operator $T_{n}$ is comonotone additive, then
$(T_{n}(f))_{n}$ converges to $f$ uniformly on the compact subsets of $X$, for
every $f\in E\cap C_{b}\left(  X\right)  $.

Notice that in both cases $(i)$ and $(ii)$ the family of testing functions can
be reduced to $1,~-\operatorname*{pr}_{1},...,~-\operatorname*{pr}_{N}$ and
$\sum_{k=1}^{N}\operatorname*{pr}_{k}^{2}$ when $K$ is included in the
positive cone of $\mathbb{R}^{N}$. Also, the convergence of $(T_{n}(f))_{n}$
to $f$ $\ $is uniform on $X$ when $f\in E$ is uniformly continuous and bounded
on $X.$
\end{theorem}

The details of this result make the objective of Section 2.

Theorem \ref{Thm2} extends not only Korovkin's original result (which
represents the particular case where $N=1,$ $K=[0,1],$ all operators $T_{n}$
are linear bounded and monotone, and the function $\operatorname*{pr}_{1}$ is
the identity of $K)$ but also the several variable version of it due to due to
Volkov \cite{V}. It encompasses also the technique of smoothing kernels, in
particular Weierstrass' argument for the Weierstrass approximation theorem:
for every bounded uniformly continuous function $f:\mathbb{R\rightarrow R},$%
\[
\left(  W_{h}f\right)  (t)=\frac{1}{h\sqrt{\pi}}\int_{-\infty}^{\infty
}f(s)e^{-\left(  s-t\right)  ^{2}/h^{2}}ds\longrightarrow f(t)
\]
uniformly on $\mathbb{R}$ as $h\rightarrow0.$

Applications of Theorem 2 in the nonlinear setting are presented in Section 3.
They are all based on Choquet's theory on integration with respect to a
capacity. Indeed, this theory, which was initiated by Choquet \cite{Ch1954},
\cite{Ch1986} in the early 1950s, represents a major source of comonotonic
additive, sublinear and monotone operators.

It is worth mentioning that nowadays Choquet's theory provides powerful tools
in decision making under risk and uncertainty, game theory, ergodic theory,
pattern recognition, interpolation theory and very recently on transport under
uncertainty. See Adams \cite{Adams},\ Denneberg \cite{Denn}, F\"{o}llmer and
Schied \cite{FS}, Wang and Klir \cite{WK}, Wang and Yan \cite{WY}, Gal and
Niculescu \cite{Gal-Nic} as well as the references therein.

For the convenience of reader we summarized in the Appendix at the end of this
paper some basic facts concerning this theory.

Some nonlinear extension of Korovkin's theorem within the framework of compact
spaces are presented in Section 4.

\section{Proof of Theorem 2}

Before to detail the proof of Theorem 2 some preliminary remarks on the
behavior of operators $T:C_{b}(X)\rightarrow C_{b}(Y)$ are necessary.

If $T$ is subadditive and monotone, then it verifies the inequality%
\begin{equation}
\left\vert T(f)-T(g)\right\vert \leq T\left(  \left\vert f-g\right\vert
\right)  \text{\quad for all }f,g. \label{f1}%
\end{equation}
Indeed, $f\leq g+\left\vert f-g\right\vert ~$yields $T(f)\leq T(g)+T\left(
\left\vert f-g\right\vert \right)  ,$ i.e., $T(f)-T(g)\leq T\left(  \left\vert
f-g\right\vert \right)  $, and interchanging the role of $f$ and $g$ we infer
that $-\left(  T(f)-T(g)\right)  \leq T\left(  \left\vert f-g\right\vert
\right)  .$

If $T$ is linear, then the property of monotonicity is equivalent to that of
positivity, whose meaning is%
\[
T(f)\geq0\text{\quad for all }f\geq0.
\]

If the operator $T$ is monotone and positively homogeneous then necessarily%
\[
T(0)=0.
\]

Every positively homogeneous and comonotonic additive operator $T$ verifies
the formula%
\begin{equation}
T(f+a\cdot1)=T(f)+aT(1)\text{\quad for all }f\text{ and all }a\in
\lbrack0,\infty); \label{f2}%
\end{equation}
indeed, $f$ is comonotonic to any constant function.

\begin{proof}
[Proof of Theorem 2.]$(i)$ In order to make more easier the handling of the
test functions we denote
\[
e_{0}=1,\text{ }e_{k}=\operatorname*{pr}\nolimits_{k}\text{ (}k=1,...N)\text{
and }e_{N+1}=\sum_{k=1}^{N}\operatorname*{pr}\nolimits_{k}^{2}.
\]

Replacing each operator $T_{n}$ by $T_{n}/T_{n}(e_{0}),$ we may assume that
$T_{n}(e_{0})=1$ for all $n.$

Let $f\in E\cap C_{b}(\Omega)$ and let $K$ be a compact subset of $X.$ Then
for every $\varepsilon>0$ there is $\tilde{\delta}>0$ such that
\[
|f(s)-f(t)|\leq\varepsilon\text{\quad for every }t\in K\text{ and }s\in
X\text{ with }\Vert s-t\Vert\leq\tilde{\delta};
\]
this can be easily proved by reductio ad absurdum.

If $\Vert s-t\Vert\geq\tilde{\delta}$, then
\[
|f(s)-f(t)|\leq\frac{2\Vert f\Vert_{\infty}}{\tilde{\delta}^{2}}\cdot\Vert
s-t\Vert^{2},
\]
so that letting $\delta=2\Vert f\Vert_{\infty}/\tilde{\delta}^{2}$ we obtain
the estimate
\begin{equation}
|f(s)-f(t)|\leq\varepsilon+\delta\cdot\Vert s-t\Vert^{2} \label{est1}%
\end{equation}
for all $t\in K$ and $s\in X.$ Since $K$ is a compact set, it can embedded
into an $N$-dimensional cube $[a,b]^{N}$ for suitable $b\geq0\geq a$ and the
estimate (\ref{est1}) yields%
\begin{multline*}
\left\vert f(s)-f(t)e_{0}\right\vert \leq\varepsilon e_{0}\\
+\delta(\varepsilon)\left[  e_{N+1}^{2}(s)+2\sum_{k=1}^{N}\left(
e_{k}(t)-a\right)  \left(  -e_{k}(s)\right)  \right. \\
\left.  -2a\sum_{k=1}^{N}e_{k}(s)+\left\Vert t\right\Vert ^{2}e_{0}(s)\right]
.
\end{multline*}
Taking into account the formula (\ref{f1}) and the fact that the operators
$T_{n}$ are subadditive and positively homogeneous, we infer that%
\begin{multline*}
\left\vert T_{n}(f)(s)-f(t)\right\vert =\left\vert T_{n}(f)(s)-T_{n}%
(f(t)e_{0})(s)\right\vert \leq T_{n}\left(  \left\vert f(s)-f(t)e_{0}%
\right\vert \right) \\
\leq\varepsilon+\delta(\varepsilon)\left[  T_{n}(e_{N+1}^{2})(s)+2\sum
_{k=1}^{N}\left(  e_{k}(t)-a\right)  T_{n}(-e_{k})(s)\right. \\
\left.  -2a\sum_{k=1}^{N}T_{n}\left(  e_{k})(s)\right)  +\left\Vert
t\right\Vert ^{2}\right]
\end{multline*}
for every $n\in\mathbb{N}$ and $s,t\in K.$ Here we used the assumption that
$f$ is nonnegative. By our hypothesis,
\[
T_{n}(e_{N+1}^{2})(s)+2\sum_{k=1}^{N}\left(  e_{k}(s)-a\right)  T_{n}%
(-e_{k})(s)-2a\sum_{k=1}^{N}T_{n}\left(  e_{k})(s)\right)  +\left\Vert
s\right\Vert ^{2}\rightarrow0
\]
uniformly on $K$ as $n\rightarrow\infty.$ Therefore
\[
\underset{n\rightarrow\infty}{\lim\sup}\left\vert T_{n}(f)(s)-f(s)\right\vert
\leq\varepsilon
\]
whence we conclude that $T_{n}(f)\rightarrow f$ uniformly on $K$ because
$\varepsilon$ was arbitrarily fixed.

$(ii)$ Suppose in addition that each operator $T_{n}$ is also comonotone
additive. According to the assertion $(i),$
\[
T_{n}(f+\Vert f\Vert e_{0})\rightarrow f+\Vert f\Vert e_{0},\text{\quad
uniformly on }K.
\]
Since a constant function is comonotone with any arbitrary function, using the
comonotone additivity of $T_{n}$ it follows that $T_{n}(f+\Vert f\Vert
e_{0})=T_{n}(f)+\Vert f\Vert\cdot T_{n}(e_{0})$. Therefore $T_{n}%
(f)\rightarrow f$ uniformly on $K$.
\end{proof}

When $K$ is included in the positive cone of $\mathbb{R}^{N}$, it can embedded
into an $N$-dimensional cube $[0,b]^{N}$ for a suitable $b>0$ and the estimate
(\ref{est1}) yields%
\begin{multline*}
\left\vert f(s)-f(t)e_{0}\right\vert \leq\varepsilon e_{0}\\
+\delta(\varepsilon)\left[  e_{N+1}^{2}(s)+2\sum_{k=1}^{N}e_{k}(t)\left(
-e_{k}(s)\right)  +\left\Vert t\right\Vert ^{2}e_{0}(s)\right]  .
\end{multline*}
Proceeding as above, we infer that%
\begin{multline*}
\left\vert T_{n}(f)(s)-f(t)\right\vert \\
\leq\varepsilon+\delta(\varepsilon)\left[  T_{n}(e_{N+1}^{2})(s)+2\sum
_{k=1}^{N}e_{k}(t)T_{n}(-e_{k})(s)+\left\Vert t\right\Vert ^{2}\right]
\end{multline*}
for every $n\in\mathbb{N}$ and $s,t\in K,$ provided that $f\geq0.$ As a
consequence, in both cases $(i)$ and $(ii)$ the family of testing functions
can be reduced to $e_{0},-e_{1},...,-e_{N}$ and $e_{N+1}.$

When dealing with functions $f\in E$ uniformly continuous and bounded on $X,$
an inspection of the argument above shows that $f$ verifies an estimate of the
form (\ref{est1}) for all $s,t\in X$, a fact that implies the convergence of
$(T_{n}(f))_{n}$ to $f$ $\ $uniformly on $X$.

\section{Applications of Theorem 2}

We will next discuss several examples of operators illustrating Theorem 2.
They are all based on Choquet's theory of integration with respect to a
capacity $\mu$, in our case the restriction of the monotone and submodular
capacity%
\[
\mu(A)=\left(  \mathcal{L}(A)\right)  ^{1/2}%
\]
to various compact subintervals of $\mathbb{R}$; here $\mathcal{L}$ denotes
the Lebesgue measure on real line. The necessary background on Choquet's
theory is provided by the Appendix at the end of this paper.

The 1-dimensional case of Theorem 2 is illustrated by the following three
families of nonlinear operators, first considered in \cite{Gal-Mjm}:

- the \emph{Bernstein-Kantorovich-Choquet }operators act on $C([0,1])$ by the
formula
\[
K_{n,\mu}(f)(x)=\sum_{k=0}^{n}\frac{(C)\int_{k/(n+1)}^{(k+1)/(n+1)}f(t)d\mu
}{\mu([k/(n+1),(k+1)/(n+1)])}\cdot{\binom{n}{k}}x^{k}(1-x)^{n-k};
\]

- the \emph{Sz\'{a}sz-Mirakjan-Kantorovich-Choquet }operators\emph{ }act on
$C([0,\infty))$ by the formula \emph{ }
\[
S_{n,\mu}(f)(x)=e^{-nx}\sum_{k=0}^{\infty}\frac{(C)\int_{k/n}^{(k+1)/n}%
f(t)d\mu}{\mu([k/n,(k+1)/n])}\cdot\frac{(nx)^{k}}{k!};
\]

- the \emph{Baskakov-Kantorovich-Choquet }operators\emph{ }act on
$C([0,\infty))$ by the formula
\[
V_{n,\mu}(f)(x)=\sum_{k=0}^{\infty}\frac{(C)\int_{k/n}^{(k+1)/n}f(t)d\mu}%
{\mu([k/n,(k+1)/n])}\cdot{\binom{n+k-1}{k}}\frac{x^{k}}{(1+x)^{n+k}}.
\]

Since the Choquet integral with respect to a submodular capacity $\mu$ is
comonotone additive, sublinear and monotone, it follows that all above
operators also have these properties.

Clearly, $K_{n,\mu}(e_{0})(x)=1$ and by Corollary 3.6 $(i)$ in \cite{Gal-Mjm}
we immediately get that $K_{n,\mu}(e_{2})(x)\rightarrow e_{2}(x)$ uniformly on
$[0,1]$. Again by Corollary 3.6 $(i),$ it follows that $K_{n,\mu}%
(1-e_{1})(x)\rightarrow1-e_{1}$, uniformly on $[0,1]$. Since $K_{n,\mu}$ is
comonotone additive,
\[
K_{n,\mu}(1-e_{1})(x)=K_{n,\mu}(e_{0})(x)+K_{n,\mu}(-e_{1})(x),
\]
which implies that $K_{n,\mu}(-e_{1})\rightarrow-e_{1}$ uniformly on $[0,1].$
Therefore the operators $K_{n,\mu}$ satisfy the hypothesis of Theorem
\ref{Thm2}, whence the conclusion%
\[
K_{n,\mu}(f)(x)\rightarrow f(x)\text{ uniformly for every }f\in C([0,1]).
\]

Similarly, one can show that the operators $S_{n,\mu}$ and $V_{n,\mu}$ satisfy
the hypothesis of Theorem 2 for $N=1$ and $X=[0,+\infty)$. In the first case,
notice that the condition $S_{n,\mu}(e_{0})=e_{0}$ is trivial. The convergence
of the sequence of functions $S_{n,\mu}(e_{2})(x)$ will be settled by
computing the integrals $\sqrt{n}\cdot(C)\int_{k/n}^{(k+1)/n}t^{2}d\mu$. We
have
\begin{multline*}
\sqrt{n}\cdot(C)\int_{k/n}^{(k+1)/n}t^{2}d\mu=\sqrt{n}\int_{0}^{\infty}%
\mu(\{t\in\lbrack k/n,(k+1)/n]:t\geq\sqrt{\alpha}\})d\alpha\\
=\sqrt{n}\int_{0}^{((k+1)/n)^{2}}\mu(\{t\in\lbrack k/n,(k+1)/n]:t\geq
\sqrt{\alpha}\})d\alpha\\
=\sqrt{n}\int_{0}^{(k/n)^{2}}\mu(\{t\in\lbrack k/n,(k+1)/n]:t\geq\sqrt{\alpha
}\})d\alpha\\
+\sqrt{n}\int_{(k/n)^{2}}^{((k+1)/n)^{2}}\mu(\{t\in\lbrack k/n,(k+1)/n]:t\geq
\sqrt{\alpha}\})d\alpha\\
=\sqrt{n}\cdot\left(  \frac{k}{n}\right)  ^{2}\cdot\frac{1}{\sqrt{n}}+\sqrt
{n}\cdot\int_{(k/n)^{2}}^{((k+1)/n)^{2}}\sqrt{(k+1)/n-\sqrt{\alpha}}d\alpha\\
=\left(  \frac{k}{n}\right)  ^{2}+\sqrt{n}\cdot\int_{0}^{1/n}\beta
^{1/2}((k+1)/n-\beta)d\beta\\
=\left(  \frac{k}{n}\right)  ^{2}+\sqrt{n}\cdot\frac{2(k+1)}{n}\cdot\frac
{2}{3}\cdot\beta^{3/2}|_{0}^{1/n}-2\sqrt{n}\cdot\frac{2}{5}\beta^{5/2}%
|_{0}^{1/n}\\
=\allowbreak\frac{1}{15n^{2}}\left(  15k^{2}+20k+8\right)  .
\end{multline*}
This immediately implies
\[
S_{n,\mu}(e_{2})(x)=S_{n}(e_{2})(x)+\frac{4}{3n}S_{n}(e_{1})(x)+\frac
{4}{3n^{2}}-\frac{4}{5n^{2}}\rightarrow e_{2}(x),
\]
uniformly on every compact subinterval $[0,a]$. Here $S_{n}$ denotes the
classical Sz\'{a}sz-Mirakjan-Kantorovich operator, associated to the Lebesgue measure.

It remains to show that $S_{n,\mu}(-e_{1})(x)\rightarrow-e_{1}(x)$, uniformly
on every compact subinterval $[0,a]$. For this goal we have to perform the
following computation:
\begin{multline*}
\sqrt{n}\cdot(C)\int_{k/n}^{(k+1)/n}(-t)d\mu=\int_{-\infty}^{0}\left\{
\mu(\{\omega\in\lbrack k/n,(k+1)/n]:-\omega\geq\alpha\})-\frac{1}{\sqrt{n}%
}\right\}  d\alpha\\
=\sqrt{n}\int_{-k/n}^{0}\left\{  \mu(\{\omega\in\lbrack k/n,(k+1)/n]:\omega
\leq-\alpha\})-\frac{1}{\sqrt{n}}\right\}  d\alpha\\
+\sqrt{n}\int_{-(k+1)/n}^{-k/n}\left\{  \mu(\{\omega\in\lbrack
k/n,(k+1)/n]:\omega\leq-\alpha\})-\frac{1}{\sqrt{n}}\right\}  d\alpha\\
=-\frac{k}{n}+\sqrt{n}\cdot\int_{-(k+1)/n}^{-k/n}\left(  \sqrt{-\alpha
-k/n}-\frac{1}{\sqrt{n}}\right)  d\alpha\\
=-\frac{k}{n}+\sqrt{n}\int_{k/n}^{(k+1)/n}\sqrt{\beta-k/n}d\beta-\frac{1}{n}\\
=-\frac{k}{n}+\sqrt{n}\int_{0}^{1/n}\beta^{1/2}d\beta-\frac{1}{n}=-\frac
{3k+1}{3n}.
\end{multline*}
Consequently
\[
S_{n,\mu}(-e_{1})(x)=S_{n}(-e_{1})(x)-\frac{1}{n}\rightarrow-x,
\]
uniformly on any compact interval $[0,a]$.

In a similar way, one can be prove that the Baskakov-Kantorovich-Choquet
operators $V_{n,\mu}$ satisfy the hypothesis of Theorem \ref{Thm2}.

The several variables framework can be illustrated by the following special
type of Bernstein-Durrmeyer-Choquet operators (see \cite{Gal-Opris} for the
general case) that act on the space of continuous functions defined on the
$N$-simplex
\[
\Delta_{N}=\{(x_{1},...,x_{N}):0\leq x_{1},...,x_{N}\leq1,\,0\leq x_{1}%
+\cdots+x_{N}\leq1\}
\]
via the formulas
\[
M_{n,\mu}(f)(\mathbf{x})=B_{n}(f)(\mathbf{x})-f(\mathbf{x})+x_{N}^{n}\left[
\frac{(C)\int_{\Delta_{N}}f(t_{1},...t_{N})t_{N}^{n}d\mu}{(C)\int_{\Delta_{N}%
}t_{N}^{n}d\mu}-f(0,...,0,1)\right]  .
\]
Here $\mathbf{x}=(x_{1},...,x_{N})$, $B_{n}(f)(\mathbf{x})$ is the
multivariate Bernstein polynomial and $\mu=\sqrt{\mathcal{L}_{N}}$, where
$\mathcal{L}_{N}$ is the $N$-dimensional Lebesgue measure. The fact that these
operators verify the hypotheses of Theorem 2 is an exercise left to the reader.

\section{\noindent The case of spaces of functions defined on compact spaces}

The alert reader has probably already noticed that the basic clue in the proof
of Theorem 2 is the estimate (\ref{est1}), characterized in \cite{N2009} (see
also \cite{N1979}) as a property of absolute continuity. This estimate occurs
in the larger context of spaces $C(M),$ where $M$ is a metric space on which
is defined a \emph{separating function}, that is, a nonnegative continuous
function $\gamma:M\times M\rightarrow\mathbb{R}$ such that%
\[
\gamma(s,t)=0\text{ implies }s=t.
\]

If $M$ is a compact subset of $\mathbb{R}^{N},$ and $f_{1},...,f_{m}\in C(M)$
is a family of functions which separates the points of $M$ (in particular this
is the case of the coordinate functions $\operatorname*{pr}_{1}%
,...,\operatorname*{pr}_{N}),$ then%
\begin{equation}
\gamma(s,t)=\sum_{k=1}^{m}\left(  f_{k}(s)-f_{k}(t)\right)  ^{2} \label{csf}%
\end{equation}
is a separating function.

\begin{lemma}
\label{LemN}\emph{(}See \emph{\cite{N2009})} If $K$ is a compact metric space,
and $\gamma:K\times K\rightarrow\mathbb{R}$ is a separating function, then any
real-valued continuous function $f$ defined on $K$ verifies an estimate of the
following form%
\[
\left\vert f(s)-f(t)\right\vert \leq\varepsilon+\delta(\varepsilon
)\gamma(s,t)\text{\quad for all }s,t\in K.
\]

\end{lemma}

The separating functions play an important role in obtaining Korovkin-type
theorems. A sample is as follows:

\begin{theorem}
Suppose that $K$ is a compact metric space and $\gamma$ is a separating
function for $M.$ If $T_{n}:C(K)\rightarrow C(K)$ $(n\in\mathbb{N})$ is a
sequence of comonotone additive, sublinear and monotone operators such that
$T_{n}(1)\rightarrow1$ uniformly and%
\begin{equation}
T_{n}(\gamma(\cdot,t))(t)\rightarrow0\text{\quad uniformly in }t, \label{cc}%
\end{equation}
then $T_{n}(f)\rightarrow f$ uniformly for each $f\in C(K).$
\end{theorem}

\noindent\ The details are similar to that used for Theorem 2, so they will be omitted.

In a similar way one can prove the following nonlinear extension of the
Korovkin type theorem (due in the linear case to Schempp \cite{Sche} and
Grossman \cite{Gro1974}):

\begin{theorem}
\emph{ }Let $X$ be a compact Hausdorff space and $\mathcal{F}$ a subset of
$C(X)$ that separates the points of $X$. If $(T_{n})_{n}$ is a sequence of
comonotonic additive, sublinear and monotone operators that map $C(X)$ into
$C(X)$ and satisfy the conditions $\lim_{n\rightarrow\infty}T_{n}(f^{k}%
)=f^{k}$ for each $f$ in $\mathcal{F}$ and $k=0,1,2,$ then%
\[
\lim_{n\rightarrow\infty}T_{n}(f)=f,
\]
for every $f$ in $C(X)$.
\end{theorem}

\section{Appendix. Some basic facts on capacities and Choquet integral}

For the convenience of the reader we will briefly recall in this section some
basic facts concerning the mathematical concept of capacity and the integral
associated to it. Full details are to be found in the books of Denneberg
\cite{Denn}, Grabisch \cite{Gr2016} and Wang and Klir \cite{WK}.

Let $(X,\mathcal{A})$\ be an arbitrarily fixed measurable space, consisting of
a nonempty abstract set $X$ and a $\sigma$-algebra ${\mathcal{A}}$ of subsets
of $X.$

\begin{definition}
\label{def1} A set function $\mu:{\mathcal{A}}\rightarrow\lbrack0,\infty)$ is
called a capacity if $\mu(\emptyset)=0$ and%
\[
\mu(A)\leq\mu(B)\text{\quad for all }A,B\in{\mathcal{A}},\text{ with }A\subset
B
\]

A capacity is called normalized if $\mu(X)=1;$
\end{definition}

An important class of normalized capacities is that of probability measures
(that is, the capacities playing the property of $\sigma$-additivity).
Probability distortions represents a major source of nonadditive capacities.
Technically, one start with a probability measure $P:\mathcal{A\rightarrow
}[0,1]$ and applies to it a distortion $u:[0,1]\rightarrow\lbrack0,1],$ that
is, a nondecreasing and continuous function such that $u(0)=0$ and
$u(1)=1;$for example, one may chose $u(t)=t^{a}$ with $\alpha>0.$The
\emph{distorted probability} $\mu=u(P)$ is a capacity with the remarkable
property of being continuous by descending sequences, that is,%
\[
\lim_{n\rightarrow\infty}\mu(A_{n})=\mu\left(
{\displaystyle\bigcap_{n=1}^{\infty}}
A_{n}\right)
\]
for every nonincreasing sequence $(A_{n})_{n}$ of sets in $\mathcal{A}.$ Upper
continuity of a capacity is a generalization of countable additivity of an
additive measure. Indeed, if $\mu$ is an additive capacity, then upper
continuity is the same with countable additivity. When the distortion $u$ is
concave (for example, when $u(t)=t^{a}$ with $0<\alpha<1),$ then $\mu$ is also
\emph{submodular} in the sense that%
\[
\mu(A\cup B)+\mu(A\cap B)\leq\mu(A)+\mu(B)\text{\quad for all }A,B\in
\mathcal{A}.
\]

Another simple technique of constructing normalized submodular capacities
$\mu$ on\ a measurable space $\left(  X,\mathcal{A}\right)  $ is by allocating
to it a probability space $\left(  Y,\mathcal{B},P\right)  $ via a map
$\rho:\mathcal{A\rightarrow B}$ such that%
\begin{gather*}
\rho(\emptyset)=\emptyset,\text{ }\rho(X)=Y\text{ and}\\
\rho\left(
{\displaystyle\bigcap\nolimits_{n=1}^{\infty}}
A_{n}\right)  =%
{\displaystyle\bigcap\nolimits_{n=1}^{\infty}}
\rho(A_{n})\text{\quad for every sequence of sets }A_{n}\in\mathcal{A}.
\end{gather*}
This allows us to define $\mu$ by the formula
\[
\mu(A)=1-P\left(  \rho(X\backslash A)\right)  .
\]
See Shafer \cite{Sha} for details.

The next concept of integrability with respect to a capacity refers to the
whole class of random variables, that is, to all functions $f:X\rightarrow
\mathbb{R}$ such that $f^{-1}(A)\in{\mathcal{A}}$ for every Borel subset $A$
of $\mathbb{R}$.

\begin{definition}
\label{def2}The Choquet integral of a random variable $f$ with respect to the
capacity $\mu$ is defined as the sum of two Riemann improper integrals,
\begin{align*}
(C)\int_{X}fd\mu &  =\int_{0}^{+\infty}\mu\left(  \{x\in X:f(x)\geq
t\}\right)  dt\\
&  +\int_{-\infty}^{0}\left[  \mu\left(  \{x\in X:f(x)\geq t\}\right)
-\mu(X)\right]  dt,
\end{align*}
Accordingly, $f$ is said to be Choquet integrable if both integrals above are finite.
\end{definition}

If $f\geq0$, then the last integral in the formula appearing in Definition
\ref{def2} is 0.

The inequality sign $\geq$ in the above two integrands can be replaced by $>;$
see \cite{WK}, Theorem 11.1,\emph{ }p. 226.

Every bounded random variable is Choquet integrable. The Choquet integral
coincides with the Lebesgue integral when the underlying set function $\mu$ is
a $\sigma$-additive measure.

The integral of a function $f:X\rightarrow\mathbb{R}$ on a set $A\in
\mathcal{A}$ is defined by the formula%
\[
(C)\int_{A}fd\mu=(C)\int_{X}fd\mu_{A}%
\]
where $\mu_{A}$ is the capacity defined by $\mu_{A}(B)=\mu(B\cap A)$ for all
$B\in\mathcal{A}.$

We next summarize some basic properties of the Choquet integral.

\begin{remark}
\label{rem1}$(a)$ If $\mu:{\mathcal{A}}\rightarrow\lbrack0,\infty)$ is a
capacity, then the associated Choquet integral is a functional on the space of
all bounded random variables such that:%
\begin{gather*}
f\geq0\text{ implies }(C)\int_{A}fd\mu\geq0\text{ \quad\emph{(}%
positivity\emph{)}}\\
f\leq g\text{ implies }\left(  C\right)  \int_{A}fd\mu\leq\left(  C\right)
\int_{A}gd\mu\text{ \quad\emph{(}monotonicity\emph{)}}\\
\left(  C\right)  \int_{A}afd\mu=a\cdot\left(  \left(  C\right)  \int_{A}%
fd\mu\right)  \text{ for }a\geq0\text{ \quad\emph{(}positive\emph{
}homogeneity\emph{)}}\\
\left(  C\right)  \int_{A}1\cdot d\mu(t)=\mu(A)\text{\quad\emph{(}%
calibration\emph{)}};
\end{gather*}
see \emph{\cite{Denn},} Proposition \emph{5.1} $(ii)$, p. \emph{64}, for a
proof of the property of positive\emph{ }homogeneity.

$(b)$ In general, the Choquet integral is not additive but, if the bounded
random variables $f$ and $g$ are comonotonic,\emph{ }then%
\[
\left(  C\right)  \int_{A}(f+g)d\mu=\left(  C\right)  \int_{A}fd\mu+\left(
\operatorname*{Ch}\right)  \int_{A}gd\mu.
\]
This is usually referred to as the property of comonotonic additivity and was
first noticed by Delacherie \cite{Del1970}. An immediate consequence is the
property of translation invariance,
\[
\left(  C\right)  \int_{A}(f+c)d\mu=\left(  C\right)  \int_{A}fd\mu+c\cdot
\mu(A)
\]
for all $c\in\mathbb{R}$ and all \hspace{0in}bounded random variables $f.$ For
details, see \emph{\cite{Denn}}, Proposition \emph{5.1, }$(vi)$, p. \emph{65.}

$(c)$ If $\mu$ is an upper continuous capacity, then the Choquet integral is
upper continuous in the sense that%
\[
\lim_{n\rightarrow\infty}\left(  \left(  C\right)  \int_{A}f_{n}d\mu\right)
=\left(  C\right)  \int_{A}fd\mu
\]
whenever $(f_{n})_{n}$ is a nonincreasing sequence of bounded random variables
that converges pointwise to the bounded variable $f.$ This is a consequence of
the Bepo Levi monotone convergence theorem from the theory of Lebesgue
integral\emph{.}

$(d)$ Suppose that $\mu$ is a submodular capacity. Then\ the associated
Choquet integral is a subadditive functional, that is,
\[
\left(  C\right)  \int_{A}(f+g)d\mu\leq\left(  C\right)  \int_{A}fd\mu+\left(
C\right)  \int_{A}gd\mu
\]
for all bounded random variables $f$ and $g.$ See \emph{\cite{Denn}},
Corollary \emph{6.4}, p. \emph{78.} and Corollary\emph{ 13.4, }p.\emph{ 161.
}It is also a submodular functional in the sense that
\[
\left(  C\right)  \int_{A}\sup\left\{  f,g\right\}  d\mu+\left(  C\right)
\int_{A}\inf\{f,g\}d\mu\leq\left(  C\right)  \int_{A}fd\mu+(C)\int_{A}gd\mu
\]
for all bounded random variables $f$ and $g.$ See \cite{Cer}, Theorem $13$
$(c)$.
\end{remark}

A characterization of Choquet integral in terms of additivity on comonotonic
functions is provided by the following analogue of the Riesz representation
theorem. See Zhou \cite{Zhou}, Theorem 1 and Lemma 3, for a simple (and more
general) argument.

\begin{theorem}
\label{thmR}Suppose that $I:C(X)\rightarrow\mathbb{R}$ is a comonotonically
additive and monotone functional with $I(1)=1$. Then it is also upper
continuous and there exists a unique upper continuous normalized capacity
$\mu:\mathcal{B}(X)\rightarrow\lbrack0,1]$ such that $I$ coincides with the
Choquet integral associated to it.

On the other hand, according to Remark \ref{rem1}, the Choquet integral
associated to any upper continuous capacity is a comonotonically additive,
monotone and upper continuous functional.
\end{theorem}

Notice that under the assumptions of Theorem \ref{thmR}, the capacity $\mu$ is
submodular if and only if the functional $I$ is submodular.

\end{document}